\documentclass[12pt,reqno]{amsart}

\newtheorem{theorem}{Theorem}[section]
\newtheorem{proposition}[theorem]{Proposition}
\newtheorem{lemma}[theorem]{Lemma}

\newtheorem*{theo}{Theorem}
\newtheorem*{theo1}{Theorem 1.6}

\theoremstyle{definition}
\newtheorem{definition}[theorem]{Definition}

\newtheorem{exam}[theorem]{Example}

\newcommand{\up}[1]{{\ensuremath{^{<{#1}>}}}}

\newcommand{\doum}[1]{\ensuremath{_{<\! \! \!<{#1}>\! \! \!>}}}
\newcommand{\main}{\ensuremath{ V_{0,i}}}
\newcommand{\sub}{\ensuremath{ V_{d,i}}}

\newcommand{\mbar}[1]{\ensuremath{\overline{M}_{\leq {#1}}}}

\newcommand{\ce}[2]{fixed\ empty\ element\ of\ ${#1}$\ for\ ${#2}$}

\begin{document}
\title{Gotzmann monomial ideals}
\author{Satoshi Murai}
\maketitle

\begin{abstract}
A Gotzmann monomial ideal of the polynomial ring is a monomial ideal which is generated in one degree and which satisfies Gotzmann's persistence theorem.
Let $R=K[x_1,\dots,x_n]$ denote the polynomial ring in $n$ variables over a field $K$ and $M^d$ the set of monomials of $R$ of degree $d$.
A subset $V\subset M^d$ is said to be a Gotzmann subset if the ideal generated by $V$ is a Gotzmann monomial ideal.
In the present paper, we find all integers $a>0$ such that every Gotzmann subset $V\subset M^d$ with $|V|=a$ is lexsegment (up to the permutation of the variables).
In addition, we classify all Gotzmann subsets of $K[x_1,x_2,x_3]$.
\end{abstract}

\section*{Introduction}
Let $K$ be an arbitrary field and $R=K[x_1,x_2,\dots,x_n]$  the polynomial ring with deg$(x_i)=1$ for $i=1,2,\dots,n$.
Let $I=\oplus_{d=0}^{\infty} I_d$ be a homogeneous ideal of $R$.
We denote the Hilbert function of $I$  by $H(I,d)$,
i.e. $H(I,d)=$dim$_K I_d$.
Let $M$ denote the set of variables $\{ x_1,x_2,\dots,x_n\}$, 
$M^d$ the set of all monomials  of degree $d$, where $M^0=\{1\}$,
and  $\overline{M_i}=M\setminus \{x_i\}$.
For a monomial $u\in R$ and a subset $V\subset M^d$,
we define
$uV=\{uv|v\in V\}$ and
$MV=\{x_iv |v\in V,\ i=1,2,\dots,n\}$.
For a finite set $V\subset M^d$,
we write $|V|$ for a number of the elements of $V$. 
Let $\gcd(V)$ denote the greatest common divisor of the monomials belonging to $V$.

Gotzmann's  persistence theorems  \cite{G} determine the growth of the Hilbert function of a homogeneous ideal for $d>\!\!>0$.
Since we are interested in Gotzmann's persistence theorem,
we recall what Gotzmann's persistence theorem is.

Let $n$ and $h$ be positive integers.
Then $h$ can be written uniquely in the form,
called the $n$\textit{th binomial representation of} $h$,

\[ h={{h(n) +n}\choose n}+{{h(n-1)+n-1}\choose n-1}+\dots+{{h(i)+i}\choose i}, \]
where $h(n)\geq h(n-1)\geq \dots\geq h(i)\geq 0,\ i\geq 1$.
See [3, Lemma 4.2.6].
\medskip

If $h={h(n)+n \choose n}+{h(n-1)+n-1\choose n-1}+\dots+{h(i)+i\choose i} $ is the $n$th binomial representation of $h$,
then we define
$$h^{ < n  >}={h(n)+n+1 \choose n}+{h(n-1)+n\choose n-1}+\dots+{h(i)+i+1\choose i} ,$$
$$h_{<\!\!\!<  n>\!\!\!>}={h(n) +n-1 \choose n - 1}+{h(n-1)+n-2\choose n-2}+\dots+{h(i) +i- 1\choose i- 1}.$$

\begin{theo}[Minimal growth of Hilbert function]
Let $I$ be a homogeneous ideal of $R=K[x_1,x_2,\dots,x_n]$.
Then one has
\begin{eqnarray}
H(I,d+1)\geq H(I,d)\up{n-1}. \label{mb}
\end{eqnarray}
\end{theo}
This theorem was proved by F. H. S. Macaulay.
We refer the reader to \cite[\S 4.2]{CM} for further informations.
\begin{theo}[\textbf{Gotzmann's Persistence Theorem {[5]}}]
Let $R=K[x_1,\dots,x_n]$ and $I$ a homogeneous ideal of $R$ generated in degree $\leq d$.
If $H(I,d+1)=H(I,d)\up{n-1}$,
then $H(I,k+1)=H(I,k)\up{n-1}$ for all $k\geq d$.
\end{theo}

A monomial ideal $I\subset R$ is called a \textit{Gotzmann monomial ideal} if $I$ is generated in one degree $d$ and if $I$ satisfies $H(I,d)\up{n-1}=H(I,d+1)$.
Instead of discussing  an ideal itself, we consider its minimal set of  monomial generators.

Let $A=(a_1,a_2,\dots,a_n)$ and $B=(b_1,b_2,\dots,b_n)$ be elements of $\mathbb{Z}_{\geq 0}^n$.
The \textit{lexicographic order} on  $\mathbb{Z}_{\geq 0}^n$ is defined by $A<B$ if the leftmost nonzero entry of $B-A$ is positive.
Moreover, the lexicographic order on monomials of the same degree is defined by
${x_1}^{a_1}{x_2}^{a_2}\dots{x_n}^{a_n}<{x_1}^{b_1}{x_2}^{b_2}\dots{x_n}^{b_n}$ if $A<B$ on $\mathbb{Z}_{\geq 0}^n$.

Let $V$ be a set of monomials of degree $d$.
\begin{itemize}
\item[(i)]
$V\subset M^d$ is called a \textit{lexsegment set} if $V$ is a set of first $|V|$ monomials in lexicographic order.
Denote the lexsegment set $V$ of $K[x_1,\dots,x_n]$ in degree $d$ with $|V|=a$ by $Lex(n,d,a)$.
\item[(ii)]
$V\subset M^d$ is called a \textit{Gotzmann set} if the ideal $I$ which is generated by $V$ satisfies $H(I,d+1) =H(I,d)\up{n-1}$,
where $I=\{0\}$ if $V=\emptyset$.
In other words, $V$ is a Gotzmann set if $\mathrm{dim}MV=\mathrm{dim}V\up{n-1}$.
\item[(iii)]
$V$ is called \textit{strongly stable} if, for any monomial $u\in V$,
one has $\frac{x_i}{x_j}u\in V$ for all $i$ and $j$ with  $i<j$ and with $x_j|u$. 
\end{itemize}

A lexsegment set is Gotzmann and strongly stable.
In general, however, a Gotzmann set is not necessarily  lexsegment.
We define $V\sim V'$ if we can obtain $V'$ from $V$ by a permutation of variables.
In other words,
there exist a permutation $\pi$ of $\{1,2,\dots,n\}$ such that $\pi(V)=V'$,
where for the permutation $\pi=(\pi(1),\dots,\pi(n))$ of $\{1,2,\dots,n\}$,
we define $\pi(x_1^{a_1}\dots x_n^{a_n})=x_{\pi(1)}^{a_1}\dots x_{\pi(n)}^{a_n}$ and $\pi(V)=\{\pi(u)|u\in V\}$.

The main result of present paper is  finding all integers $ a>0$ such that
every Gotzmann set $V$  of degree $d$ with $|V|=a$ and with $\gcd(V)$ satisfies $V\sim Lex(n,d,a)$.

\begin{theo1}
Let $R=K[x_1,\dots,x_n]$ be the polynomial ring in $n$ variables
and $a=\sum_{j=p}^{n-1}{{a(j)+j}\choose j}$  the $(n-1)th$ binomial representation of $a>0$.
Then the following conditions are equivalent:
\begin{itemize}
\item[(i)]$a(n-1)=a(p)$;
\item[(ii)]For every Gotzmann set $V$ with $|V|=a$ and $\gcd(V)=1$, one has $V\sim Lex(n,d,a)$, where $V\subset M^d$;
\item[(iii)]For every Gotzmann set $V$ with $|V|=a$ and $\gcd(V)=1$, one has $V\sim V'$ for some strongly stable set $V'$ consisting of  monomials of $R$.
\end{itemize}
\end{theo1}
\medskip

In Theorem \ref{thm1}, we describe nothing about the degree $d$.
But, in Lemma \ref{degree},
we will prove that the degree $d$ of Gotzmann set $V$ with $|V|=a$ and $\gcd(V)=1$ is determined automatically.

One Related work of Gotzmann's theorems
is  Gotzmann theorems for exterior algebra which are done by A. Aramova, J. Herzog and T. Hibi  \cite{AH}.
Gotzmann set of exterior algebra is called a \textit{squarefree Gotzmann set}.
Similar consequence of Theorem \ref{thm1} was done by Z. Furedi and J. R. Griggs \cite{FG}.
They determined all integers $a>0$ such that every squarefree Gotzmann set $V$ with $|V|=a$ is unique up to the permutation of variables.

We also classify all Gotzmann sets of $K[x_1,x_2,x_3]$ in Proposition \ref{prop1}.

\section{Proof of Theorem \ref{thm1}}
Let $K$ be an arbitrary field and $R=K[x_1,\dots,x_n]$ the polynomial ring in $n$ variables over $K$.

Let $h$ be a positive integer and $h=\sum_{j=i}^n {h(j) +j \choose j}$ the $n$th binomial representation of $h$.
Let $\alpha =$max$\{0, $max$\{\alpha\in \mathbb{Z}|h-{{\alpha + n}\choose n}>0\}\}$.
We denote $h-{{\alpha +n}\choose n}$ by $\bar{h}^{(n)}$, in other words,
\begin{itemize}
\item[(i)] if $h=1$, then $\bar{h}^{(n)}=0$,
\item[(ii)] if $h>1$ and $i=n$, then $\bar{h}^{(n)}={{h(n)+n-1} \choose{n-1}}$,
\item[(iii)] if $h>1$ and $i<n$, then $\bar{h}^{(n)}=\sum_{j=i}^{n-1} {{h(j)+j}\choose j}$.
\end{itemize}

\begin{lemma}\label{gcd}
Let $V$ be a  set of monomials of same degree, $u=\gcd(V)$ and $uV'=V$.
Then $V$ is a Gotzmann set if and only if $V'$ is a Gotzmann set.
\end{lemma}
\begin{proof}
By constraction, we have $|V|=|V'|$ and $|MV|=|MV'|$.
Thus these conditions are equivalent.
\end{proof}
We can make a Gotzmann set which is not a lexsegment set by multiplying a lexsegment set by some monomial.
For example, if $V$ is a lexsegment set,
then $x_1x_2V$ is not a lexsegment set, but a Gotzmann set.
But this is essentially the same as a lexsegment set.
Therefore we often assume $\gcd(V)=1$.

Let $V$ be a set of monomials of degree $d$ and  $u=\gcd(V)$.
If $|V|>1$, we define $\main=\{v\in M^d|\ x_iu\ \mathrm{divides} \ v\} $  and $\sub =V \setminus \main$  for $i=1,2,\dots,n$.
If $ |V|=1$, then  we define \main $=V$ and $\sub=\emptyset$.
Note that if $|V|>1$ then $\sub\ne\emptyset$.

First, we need some Lemmas in \cite{S} to prove the main theorem.

\begin{lemma}[{\cite[Lemma 1.5]{S}}]\label{plus}
Let $a$, $b$ and $n$ be positive integers. One has
$$a\up{n}+b\up{n}>(a+b)\up{n}.$$
\end{lemma}

\begin{lemma}[{\cite[Lemma 2.2]{S}}]\label{hanni}
Let $V$ be a Gotzmann set of monomials of degree d.
Then, for $i =1,2,\dots,n$, we have \begin{eqnarray}
\overline {|V|}^{(n-1)}\leq  |\sub| \leq |V| \doum{n-1}.\label{A-1}
\end{eqnarray}
\end{lemma}

\begin{lemma}[{\cite[Lemma 2.3]{S}}]\label{divide}
Let $V$ be a Gotzmann set of monomials  of degree $d$ with $\gcd(V)=1$ and with $V\ne M^d$.
Then there exist $i$ which satisfies followings:
\begin{itemize}
\item[(i)] \sub \ is a Gotzmann set of $K[x_1,\dots,x_{i-1},x_{i+1},\dots,x_n]$,
\main \  is a Gotzmann set of $K[x_1,\dots,x_n]$ and $| \sub|<|V| \doum{n-1}$;
\item[(ii)]   $x_i\sub \subset \overline{M_i}\main$.
\end{itemize}
\end{lemma}

\begin{lemma}\label{degree}
Let V be a Gotzmann set of monomials of degree d.
If $\gcd(V)=1$  and ${{\alpha+n-1}\choose n-1} <|V|\leq  {{\alpha+1+n-1}\choose n-1}$,
then we have $d=\alpha+1$ and  $\gcd(\main)=x_i$ if $|V|>1$. 
\end{lemma}
\begin{proof}
We use induction on $|V|$.
In  case of $|V|=1$, we have $V=\{1\}=M^0$.
Thus we may assume $|V|>1$ and $n>1$.
By Lemma \ref{divide} we may take \main \ as a  Gotzmann set.
Since $V=\main \cup \sub$ and $\sub\ne \emptyset$, we can use induction.

Let $u=\gcd(\main)$ and  $\main=uV'$.
Lemma \ref{hanni} says \[
{{\alpha-1+n-1}\choose n-1} <|\main|=  |V|- |\sub| \leq  {{\alpha+n-1}\choose n-1}.
\]
Thus, by induction, we have $V'\subset M^{\alpha}$.
On the other hand, we have $M\main\supset x_i\sub$ by Lemma \ref{divide}.
Thus for any $v\in \sub,$ we have $x_iv\in M\main=uMV'$.
Thus $u$ divides $x_iv$.
Definition of $\main$ says $x_i$ divides $u$.
Thus $u/x_i$ divides $v$ and all elements of $V$ can be divided by $u/x_i$.
Since $\gcd(V)=1$, we have $u=x_i$.
Thus $\main=x_iV'\subset M^{\alpha+1}$.
Then we have $d=\alpha+1$.
\end{proof}

Lemma \ref{degree} says that
if we fixed $|V|$ and assume $\gcd(V)=1$ then the degree of elements of $V$ is automatically determined.
Furthermore,
if $V$ is a Gotzmann set with $|V|={{\alpha+n}\choose n-1}$ and $\gcd(V)=1$, then $V=M^{\alpha+1}$.

\begin{theorem}\label{thm1}
Let $R=K[x_1,\dots,x_n]$ be the polynomial ring in $n$ variables
and $a=\sum_{j=p}^{n-1}{{a(j)+j}\choose j}$ the $(n-1)th$ binomial representation of $a>0$.
Then the following conditions are equivalent:
\begin{itemize}
\item[(i)]$a(n-1)=a(p)$;
\item[(ii)]For every Gotzmann set $V$ with $|V|=a$ and $\gcd(V)=1$, one has $V\sim Lex(n,d,a)$, where $V\subset M^d$;
\item[(iii)]For every Gotzmann set $V$ with $|V|=a$ and $\gcd(V)=1$, one has $V\sim V'$ for some strongly stable set $V'$ consisting of  monomials of $R$.
\end{itemize}
\end{theorem}

\begin{proof}
((i)$\Rightarrow$ (ii))
Let $a=a(n-1)=a(n-2)=\dots=a(p)$.
We use induction on $|V|$.
If $|V|=1$, then $V=\{1\}$ since $\gcd(V)=1$.
Thus $V$ is a lexsegment set.
If $n-1=p$, then $V=M^d$ by Lemma \ref{degree}.
Thus we may assume $p<n-1$ and  $|V|>1$.
By Lemma \ref{divide}, we may assume $|\sub|<|V|\doum{n-1}$ and \sub \ is a Gotzmann set.
Thus, by Lemma \ref{hanni}, we have the form

\begin{eqnarray*}
|\sub|=\sum_{j=p}^{n-2} {{a+j}\choose j} +b\ \ \mbox{and} \ \  |\main|=\sum_{j=p+1}^{n-1} {a-1+j \choose j} + {a+p-1\choose p} +c \\
\end{eqnarray*}\vspace{-10pt}
with $0\leq b < { a+p-1 \choose p-1}$ and $0  <c \leq {a+p-1\choose p-1}$.

Since Lemma \ref{divide} (ii) says $M\main\supset x_i \sub$,
$MV=M\main\cup \overline{M_i}\sub$ is a disjoint union.
Then by Lemmas  \ref{plus} and \ref{divide}, we have
\begin{eqnarray*}
|MV|&=&  |M\main| +|\overline{M_i} \sub| \\
&=&\bigg\{\sum_{j=p}^{n-2} {{a+j}\choose j}\bigg\}^{[+1]} + \bigg\{ \sum_{j=p}^{n-1} {a-1+j \choose j}\bigg\}^{[+1]} +b \up {p-1} + c \up {p-1}\\
&\leq&  \sum_{j=p+1}^{n-1} {a+j \choose j} +{a+p\choose p} + \{b +c\}^\up{p-1} = |V| \up{n-1}.
\end{eqnarray*} 
Lemma \ref{plus} says this is equal if and only if $b=0$ or $c=0$.
Thus $b=0$.

Then we have $|\main|={{a+n-1}\choose n-1}$.
Thus \main $=x_iM^{d-1}$.
Moreover, $|\sub|=\sum_{j=p}^{n-2}{{a(j)+j}\choose j}$ and \sub \ is a Gotzmann set of $n-1$ variables.
Thus by induction, \sub \  is a lexsegment set by proper permutation of variables.
We may assume \sub  \ is a lexsegment set of $K[x_2,\dots,x_n]$ and $i=1$.
Since $V_{d,1}$ is a lexsegment set of $K[x_2,\dots,x_n]$,
$V=x_1V_{0,1}\cup V_{d,1}=x_1M^{d-1}\cup V_{d,1}$ is a lexsegment set.
\medskip

\noindent ((ii)$\Rightarrow$(iii))
Since lexsegment sets are strongly stable,
(ii) $\Rightarrow$ (iii) is obvious.
\medskip

\noindent((iii)$\Rightarrow$(i))
In case of $a(n-1)>a(p)$,
we will construct a Gotzmann set that is not strongly stable.
By assumption, we have $a(n-1)>a(p)$.
Thus there exists $k$ such that $a(k+1)>a(k)$.
Let $V_{n-1}=x_nM^{a(n-1)}=u_{n-1}M^{a(n-1)}$.
Denote $\{x_1,x_2,\dots,x_j\}$ by \mbar{j} .
Inductively we define $V_j$ as follows:\begin{itemize}
\item If $j\ne k$, then we define  $V_j=u_j\mbar{j+1}^{a(j)}$, where $u_j=u_{j+1}\frac{{x_{j+1}}^{1+a(j+1)-a(j)}}{x_{j+2}}$;\\
\item If $j=k$, then we define $V_k=u_k\mbar{k+1}^{a(k)}$, where $u_k=u_{k+1}x_1 \frac{{x_{k+1}}^{a(k+1)-a(k)}}{x_{k+2}}$.\\
\end{itemize}
Let $V=\bigcup_{j=p}^{n-1} V_j$.
If $i>j$, then we have $V_j\cap V_i=\emptyset$
since $V_j$ has no element which can be divided by $u_i$,
Thus $V=\bigcup_{j=p}^{n-1} V_j$ is a disjoint union,
therefore, $|V|=\sum_{j=p}^{n-1} |V_j|=\sum_{j=p}^{n-1}{{a(j)+j}\choose j}=a$.
Moreover, since $u_{j+1}|x_{j+2}u_j$, we have $u_i|x_{i+1}x_{i}\cdots x_{j+2}u_j$ for $i>j$.
Since $u_i\in K[x_1,x_{i+1},x_{i+2},\dots,x_n]$, for $i>j$ we have
$$x_{i+1}V_j
\subset \frac{x_{i+1}}{x_i} \frac{x_i}{x_{i-1}} \cdots \frac{x_{j+3}}{x_{j+2}} x_{j+2}u_j\overline{M}_{\leq i+1}^{a(j)} \subset u_i\overline{M}_{\leq i+1}^{a(i)+1}
\subset\overline{M}_{\leq i+1}V_i.$$
Hence, we have
$MV=\bigcup_{j=p}^{n-1} MV_j=\bigcup_{j=p}^{n-1} \overline{M}_{\leq j+1}V_j$.
This union is also disjoint.
Thus we have $|MV|=\sum_{j=p}^{n-1} {{a(j)+1+j}\choose j}=a\up{n-1}$
and $V$ is a Gotzmann set. 

Next, we will prove $V$ is not strongly stable.
Let $u'={u_{k+1}}/{x_{k+2}}\in K[x_{k+2},\dots,x_{n}]$.
Since $a(k+1)>a(k)$, we have deg$(u')=\deg(u_k)-1-(a(k+1)-a(k))\leq d-2$.
Let $d_0=$deg$(u')$.
Then we will prove  $u'{x_1}^{d-d_0}$ and $u'{x_{k+1}}^{d-d_0}$ do not belong to $V$.
What we have to prove is $u'x_1^{d-d_0}\notin V_j$ and $u'{x_{k+1}}^{d-d_0}\notin V_j$ for all $j$.
\begin{itemize}
\item[(i)] For $j=k$, since $V_k=x_1x_{k+1}^{a(k+1)-a(k)}u'\mbar{k+1}^{a(k)}$,
we have $u'x_1^{d-d_0}\notin V_k$ and $u'{x_{k+1}}^{d-d_0}\notin V_k$.
\item[(ii)] For any $j<k$, $x_{j+1}$ divides  $u\in V_j$.
Since $u'\in K[x_{k+1},\dots,x_n]$, it follows that $u'x_1^{d-d_0}\notin V_j$ and $u'{x_{k+1}}^{d-d_0}\notin V_j$.
\item[(iii)] For any $j>k$, $u_j$ does not divide $u'$.
Since $u_j\in K[x_{j+1},\dots,x_{n}]$, $u_j$ does not divide $u'{x_1}^{d-d_0}$ and $u'{x_{k+1}}^{d-d_0}$.
Thus $u'x_1^{d-d_0}\notin V_j$ and $u'{x_{k+1}}^{d-d_0}\notin V_j$.
\end{itemize}
However, if $V$ is strongly stable,
then either $u'x_1^{d-d_0}\in V$ or $u'x_{k+1}^{d-d_0}\in V$ must be satisfied  since  $u'x_1{x_{k+1}}^{d-d_0-1}\in V_k\subset V$.
Thus $V$ is not strongly stable. 
\end{proof}
\medskip

\begin{definition}
Let $V$ be a Gotzmann set of monomials and  $|V|=\sum_{j=p}^{n-1} {{a(j)+j}\choose j}$  the $(n-1)$th binomial representation.
By Theorem \ref{thm1}, if $a(p)=a(n-1)$ then $V$ must be a lexsegment set.
We call $|V|$  a $n$\textit{th lexnumber}, or simply a \textit{lexnumber}  if $a(p)=a(n-1)$.
\end{definition}
\begin{exam}
Here are some lexnumbers for $n=3,4,5$.
\begin{itemize}
\item[]
\begin{itemize}
\item[$n=3$:]$1,2,3,5,6,9,10,14,15,20,21,27,28,35,36,44,45,54,55,65,66,$\\
$77,78,90,91,104,105,\dots $
\item[$n=4$:]$1,2,3,4,7,9,10,16,19,20,30,34,35,50,55,56,77,83,84,112,\dots,$
\item[$n=5$:]$1,2,3,4,5,9,12,14,15,25,31,34,35,55,65,104,105,\dots,$
\end{itemize}
\end{itemize}
If we fixed $d$, then there are only $\{d(n-1)+1\}$ lexnumbers,
since there are $(n-1)$ lexnumbers between ${{t+n-1}\choose n-1}$ and ${{t+n}\choose n-1}$.
\end{exam}
\medskip

\section{Gotzmann sets of three variables}

In this section we consider Gotzmann sets of a few variables.
If $n=1$, then all sets $V$ are Gotzmann sets.
If $n=2$, we can easily show $V$ is a Gotzmann sets if and only if $V=\emptyset$ or $V=M^d$,
when we assume $\gcd(V)=1$.
We consider the case $n=3$ in Proposition \ref{prop1}.

We define a map $\pi_i$ : $\bigoplus_{d=0}^\infty M^d \rightarrow \mathbb{Z}_{\geq 0}^{\ n-1}$ by setting
\[\pi_i(x_1^{a_1}\dots x_n^{a_n})=(a_1,\dots,a_{i-1},a_{i+1},\dots,a_n).\]
It follows that $\pi_i |_{M^d}$ is injective.

Let $V$ be a set of monomials of degree $d$ and let $u=x_1^{a_1}x_2^{a_2}\dots x_n^{a_n}$  a monomial of degree $d$.
We say that a monomial $v=x_1^{b_1}x_2^{b_2}\dots x_n^{b_n}$ with degree $d$ is \textit{under} $u$ \textit{for} $i$ 
if for any $j\ne i$, $b_j\leq a_j$.
We call $u$  a \textit{fixed empty element of} $V$ \textit{for} $i$  if $u\notin V$ and  any monomial which is under $u$ for $i$ does not belongs to $V$.

Note that $u$ is a \ce{V}{i} if and only if  $\pi_i(u)\notin\pi_i(M^tV)$ for $t\geq 0$.
Furthermore, if $u$ is a \ce{V}{i} , then any monomial $v$ which is under $u$ for $i$ is a \ce{V}{i}.
\medskip

\begin{proposition}\label{prop1}
Let $V$ be a set of monomials of degree $d$ and  $u=\gcd(V)$.
If $V$ is a Gotzmann set,
then any monomial $v\notin V$  is a \ce{V}{ \mathit{some} \ i} and $|V| >{{d-\deg(u)-1+n-1}\choose n-1}$.

Especially, if $n=3$, then these conditions are equivalent.
\end{proposition}
\begin{proof}
Let $u=\gcd(V)$ and  $ V=uV'$.
Then $\gcd(V')=1$.
By Lemma \ref{degree}, we have $|V|=|V'|>{{d-\deg(u)-1+n-1}\choose n-1}$.
First, we will prove that we may assume $\gcd(V)=1$.

Let $u=x_1^{a_1}x_2^{a_2}\dots x_n^{a_n}$.
For any monomial $x_1^{b_1}x_2^{b_2}\dots x_n^{b_n}$ with degree $d$, 
if $u$ does not divide $ x_1^{b_1}x_2^{b_2}\dots x_n^{b_n} $ then there exists  $i$ such that $b_i<a_i$.
Furthermore, for any $j\ne i$,  $x_1^{b_1}x_2^{b_2}\dots x_n^{b_n}$ is \ce{V}{\mathrm{some} \ j}.
Hence all monomials which are not divided by $u$ are fixed empty elements.
Now, if $v$ is a \ce{V'}{i},
then $uv$ is a \ce{V}{i},
because all monomials which cannot be divided by $u$ are not belong to $V$.
Hence we may  only consider $V'$.
Thus we may assume $\gcd(V)=1$.

We use induction on $|V|$.
If $|V|=1$, then $V= \{ 1 \}$.
Thus, in this case,  the conditions are satisfied.
Hence we may assume $|V|>1$.
By Lemma \ref{divide}, there exists  $i$ such that \main \ and \sub \  are Gotzmann sets and  $\overline{M_i} \main \supset x_i \sub$.
We consider two cases for $w\not\in V$.

Let $w=x_1^{a_1}x_2^{a_2}\dots x_n^{a_n}$ be a monomial of degree $d$.
We will prove  that if $w\notin V$ then  $w$ is a \ce{V}{\mathrm{some} \ j}.
\medskip

\noindent \textbf{[Case I]}
If $x_i$ divide $w$, then $w \in \main$.
By induction, there exist  $j$ such that $w$ is  a \ce{\main}{j}.
Thus for any $v\ne w$ being under $w$ for $j$, we have $v\notin \main$.
Hence what we have to prove is $v\notin \sub$.

If $j=i$ then $x_i|v$, thus $v\notin \sub$.

If $j\ne i$, we may assume $x_i\not | v$ since if $x_i | v$ then $v\notin \sub$.
Let $v=x_1^{b_1}x_2^{b_2} \dots x_n^{b_n}\ne w$.
Since $x_i|w$ and $x_i \not |v$, we have $b_i+1\leq a_i$.
Thus  for any  $x_k|v$, $\frac{x_i}{x_k} v$ is under $u$ for $j$.
Hence we have $\frac{x_i}{x_k} v\not\in\main$ and $x_iv\notin M\main$.
Since $\overline{M_i}\main\supset x_i\sub$,
we have $v\notin \sub$.
\medskip

\noindent \textbf{[Case II]} If $x_i$ does not divide $w$, then by induction there exists  $j\ne i$ such that $w$ is a \ce{\sub}{j}.
Since $j\ne i$, for any $v$ being under $w$ for $i$, $v$ cannot be divided by $x_i$.
Thus we have $v \notin \main \cup \sub=V$.
\medskip

Next, in case of $n=3$, we will prove these conditions are equivalent.
We may assume $\gcd(V)=1$.
Let $u=x_1^{a_1}x_2^{a_2}x_3^{a_3}$ be a \ce{V}{i}.
We consider the case of $i=1$.
Since $x_1^{a_1+1}x_2^{a_2-1}x_3^{a_3}$ and $x_1^{a_1+1}x_2^{a_2}x_3^{a_3-1}$ are under $u$ for $1$ and
$u\notin V$,
we have $x_1u=x_1^{a_1+1}x_2^{a_2}x_3^{a_3}\notin MV$.
By the same way, if $u$ is a \ce{V}{i} then  $x_i u$ is a \ce{MV}{i}. 
For any $u'\ne u$ which is a \ce{V}{j},
we will prove $x_iu\ne x_ju'$.

If $i=j$, then we have $x_iu\ne x_ju'$ since $u\ne u'$.
If $i\ne j$, then  monomials $u=x_i^{a_i}x_j^{a_j}x_k^{a_k},x_i^{a_i+1}x_j^{a_j-1}x_k^{a_k},\dots,x_i^{a_i+a_j}x_j^{0}x_k^{a_k}$
are under $u$ for $i$.
Moreover, monomials $x_i^{a_i-1}x_j^{a_j+1}x_k^{a_k},x_i^{a_i-2}x_j^{a_j+2}x_k^{a_k},\dots,x_j^{a_j+a_i}x_k^{a_k},
x_j^{a_j+a_i+1}x_k^{a_k-1},\dots,x_j^{d}$ are under $u'$ for $j$.
Hence  we can take $\{ (a_j+1)+(d-a_j) \}$ monomials which do not belong to $V$.
Thus we have
$|V| \leq |M^d|-(d+1)={{d+1}\choose2}.$
But  assumption says $|V|>{{d+1}\choose 2}$, this is a contradiction.
Thus if $u$ is a \ce{V}{i} and $u'\ne u$ is a \ce{V}{j},
then we have $x_iu\notin MV$, $x_ju'\notin MV$ and $x_iu\ne x_ju'$.

Hence if $V$ has $l$ fixed empty elements,
then $MV$ has at least $l$ fixed empty elements.
Thus  we have 
$$|MV| \leq {{d+3}\choose 2}-l={{d+2}\choose 2}+{{d+2-l}\choose 1}.$$
Moreover, by the minimal growth of Hilbert function (\ref{mb}),
we have
$$|MV|\geq |V|\up{2}=\bigg\{{{d+2}\choose 2}-l\bigg\}\up{2}=\bigg\{ {{d+1}\choose 2}+{{d+1-l}\choose 1}\bigg\}\up{2}.$$
We have $|MV|=|V|\up{2}$.
Thus $V$ is a Gotzmann set. 
\end{proof}
\begin{exam}
To understand the meaning of Proposition \ref{prop1},
drawing a picture of monomials is useful.
(Similar idea could be found  in  \cite{Gr}.)
In the picture below, we show all monomials in $K[x_1,x_2,x_3]$ with degree $4$.
The monomial $x_1^4$ is in the lower left corner, $x_3^4$ is in the lower right corner,
and $x_2^4$ is at the top.
The black dots denote monomials in $V$ and the empty circles denote monomials which are missing.
For example,
figure (1) means $x_1^4,\ x_1^3x_2,\ x_1^2x_2^2$ and $x_1 x_2^3$ are missing. 
In the picture below,
we classify all Gotzmann set $V$ in $K[x_1,x_2,x_3]$ with $\gcd(V)=1$ and $|V|={4+2 \choose 2} - 4=11$
up to permutation.
\end{exam}

\hspace{15pt}(1) \hspace{86pt} (2) \hspace{86pt} (3) \hspace{86pt} (4)

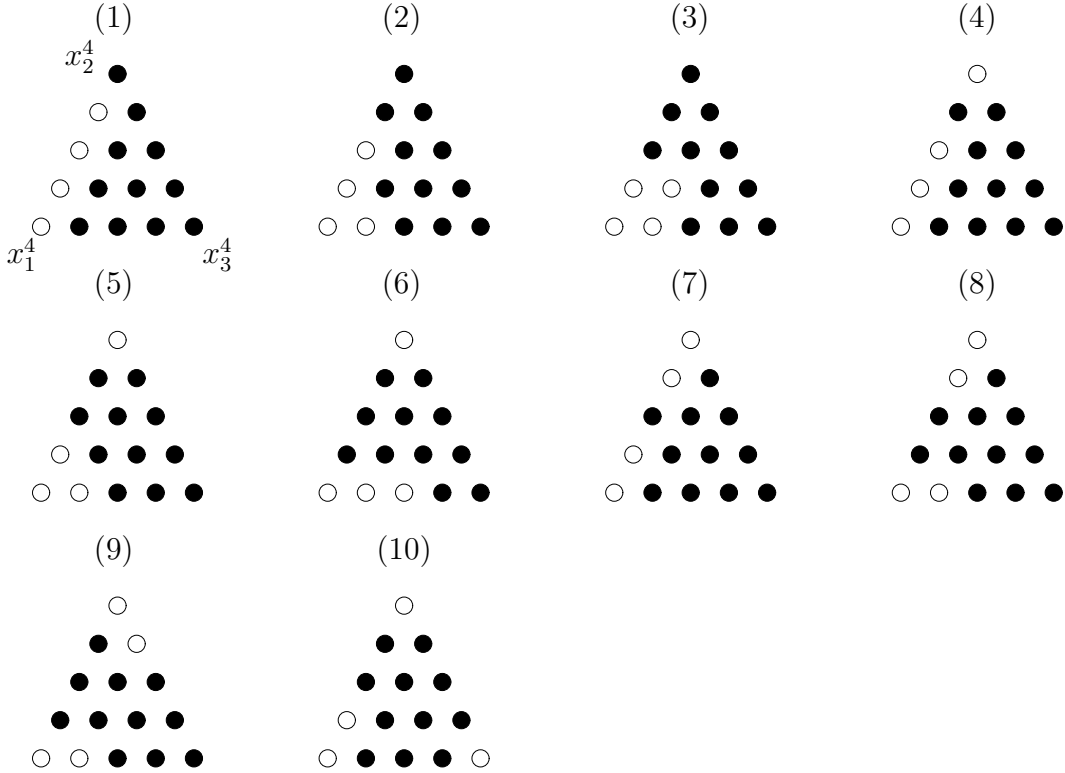
\begin{picture}(36.47,22.63)(12.25,-30.63)
\put(22,-15){\makebox(0,0){$x_2^4$}}%
\put(0,-90.){\makebox(0,0){$x_1^4$}}%
\put(74,-90){\makebox(0,0){$x_3^4$}}%
\special{sh 1.000}%
\special{ar 500 300 46 46  0.0000000 6.2831853}%

\special{ar 400 500 46 46  0.0000000 6.2831853}%
\special{sh 1.000}%
\special{ar 600 500 46 46  0.0000000 6.2831853}%

\special{ar 300 700 46 46  0.0000000 6.2831853}%
\special{sh 1.000}%
\special{ar 500 700 46 46  0.0000000 6.2831853}%
\special{sh 1.000}%
\special{ar 700 700 46 46  0.0000000 6.2831853}%

\special{ar 200 900 46 46  0.0000000 6.2831853}%
\special{sh 1.000}%
\special{ar 400 900 46 46  0.0000000 6.2831853}%
\special{sh 1.000}%
\special{ar 600 900 46 46  0.0000000 6.2831853}%
\special{sh 1.000}%
\special{ar 800 900 46 46  0.0000000 6.2831853}%

\special{ar 100 1100 46 46  0.0000000 6.2831853}%
\special{sh 1.000}%
\special{ar 300 1100 46 46  0.0000000 6.2831853}%
\special{sh 1.000}%
\special{ar 500 1100 46 46  0.0000000 6.2831853}%
\special{sh 1.000}%
\special{ar 700 1100 46 46  0.0000000 6.2831853}%
\special{sh 1.000}%
\special{ar 900 1100 46 46  0.0000000 6.2831853}%

\special{sh 1.000}%
\special{ar 2000 300 46 46  0.0000000 6.2831853}%

\special{sh 1.000}%
\special{ar 1900 500 46 46  0.0000000 6.2831853}%
\special{sh 1.000}%
\special{ar 2100 500 46 46  0.0000000 6.2831853}%

\special{ar  1800 700 46 46  0.0000000 6.2831853}%
\special{sh 1.000}%
\special{ar 2000 700 46 46  0.0000000 6.2831853}%
\special{sh 1.000}%
\special{ar 2200 700 46 46  0.0000000 6.2831853}%

\special{ar 1700 900 46 46  0.0000000 6.2831853}%
\special{sh 1.000}%
\special{ar 1900 900 46 46  0.0000000 6.2831853}%
\special{sh 1.000}%
\special{ar 2100 900 46 46  0.0000000 6.2831853}%
\special{sh 1.000}%
\special{ar 2300 900 46 46  0.0000000 6.2831853}%

\special{ar 1600 1100 46 46  0.0000000 6.2831853}%

\special{ar 1800 1100 46 46  0.0000000 6.2831853}%
\special{sh 1.000}%
\special{ar 2000 1100 46 46  0.0000000 6.2831853}%
\special{sh 1.000}%
\special{ar 2200 1100 46 46  0.0000000 6.2831853}%
\special{sh 1.000}%
\special{ar 2400 1100 46 46  0.0000000 6.2831853}%

\special{sh 1.000}%
\special{ar 3500 300 46 46  0.0000000 6.2831853}%

\special{sh 1.000}%
\special{ar 3400 500 46 46  0.0000000 6.2831853}%
\special{sh 1.000}%
\special{ar 3600 500 46 46  0.0000000 6.2831853}%

\special{sh 1.000}%
\special{ar  3300 700 46 46  0.0000000 6.2831853}%
\special{sh 1.000}%
\special{ar 3500 700 46 46  0.0000000 6.2831853}%
\special{sh 1.000}%
\special{ar 3700 700 46 46  0.0000000 6.2831853}%

\special{ar 3200 900 46 46  0.0000000 6.2831853}%

\special{ar 3400 900 46 46  0.0000000 6.2831853}%
\special{sh 1.000}%
\special{ar 3600 900 46 46  0.0000000 6.2831853}%
\special{sh 1.000}%
\special{ar 3800 900 46 46  0.0000000 6.2831853}%

\special{ar 3100 1100 46 46  0.0000000 6.2831853}%

\special{ar 3300 1100 46 46  0.0000000 6.2831853}%
\special{sh 1.000}%
\special{ar 3500 1100 46 46  0.0000000 6.2831853}%
\special{sh 1.000}%
\special{ar 3700 1100 46 46  0.0000000 6.2831853}%
\special{sh 1.000}%
\special{ar 3900 1100 46 46  0.0000000 6.2831853}%


\special{ar 5000 300 46 46  0.0000000 6.2831853}%

\special{sh 1.000}%
\special{ar 4900 500 46 46  0.0000000 6.2831853}%
\special{sh 1.000}%
\special{ar 5100 500 46 46  0.0000000 6.2831853}%

\special{ar  4800 700 46 46  0.0000000 6.2831853}%
\special{sh 1.000}%
\special{ar 5000 700 46 46  0.0000000 6.2831853}%
\special{sh 1.000}%
\special{ar 5200 700 46 46  0.0000000 6.2831853}%

\special{ar 4700 900 46 46  0.0000000 6.2831853}%
\special{sh 1.000}%
\special{ar 4900 900 46 46  0.0000000 6.2831853}%
\special{sh 1.000}%
\special{ar 5100 900 46 46  0.0000000 6.2831853}%
\special{sh 1.000}%
\special{ar 5300 900 46 46  0.0000000 6.2831853}%

\special{ar 4600 1100 46 46  0.0000000 6.2831853}%
\special{sh 1.000}%
\special{ar 4800 1100 46 46  0.0000000 6.2831853}%
\special{sh 1.000}%
\special{ar 5000 1100 46 46  0.0000000 6.2831853}%
\special{sh 1.000}%
\special{ar 5200 1100 46 46  0.0000000 6.2831853}%
\special{sh 1.000}%
\special{ar 5400 1100 46 46  0.0000000 6.2831853}%
\end{picture}

\vspace{60pt}
\hspace{15pt}(5) \hspace{86pt} (6) \hspace{86pt} (7) \hspace{86pt} (8)

\begin{picture}(36.47,22.63)(12.25,-30.63)

\special{ar 500 300 46 46  0.0000000 6.2831853}%
\special{sh 1.000}%
\special{ar 400 500 46 46  0.0000000 6.2831853}%
\special{sh 1.000}%
\special{ar 600 500 46 46  0.0000000 6.2831853}%
\special{sh 1.000}%
\special{ar 300 700 46 46  0.0000000 6.2831853}%
\special{sh 1.000}%
\special{ar 500 700 46 46  0.0000000 6.2831853}%
\special{sh 1.000}%
\special{ar 700 700 46 46  0.0000000 6.2831853}%

\special{ar 200 900 46 46  0.0000000 6.2831853}%
\special{sh 1.000}%
\special{ar 400 900 46 46  0.0000000 6.2831853}%
\special{sh 1.000}%
\special{ar 600 900 46 46  0.0000000 6.2831853}%
\special{sh 1.000}%
\special{ar 800 900 46 46  0.0000000 6.2831853}%

\special{ar 100 1100 46 46  0.0000000 6.2831853}%

\special{ar 300 1100 46 46  0.0000000 6.2831853}%
\special{sh 1.000}%
\special{ar 500 1100 46 46  0.0000000 6.2831853}%
\special{sh 1.000}%
\special{ar 700 1100 46 46  0.0000000 6.2831853}%
\special{sh 1.000}%
\special{ar 900 1100 46 46  0.0000000 6.2831853}%


\special{ar 2000 300 46 46  0.0000000 6.2831853}%

\special{sh 1.000}%
\special{ar 1900 500 46 46  0.0000000 6.2831853}%
\special{sh 1.000}%
\special{ar 2100 500 46 46  0.0000000 6.2831853}%
\special{sh 1.000}%
\special{ar  1800 700 46 46  0.0000000 6.2831853}%
\special{sh 1.000}%
\special{ar 2000 700 46 46  0.0000000 6.2831853}%
\special{sh 1.000}%
\special{ar 2200 700 46 46  0.0000000 6.2831853}%

\special{sh 1.000}%
\special{ar 1700 900 46 46  0.0000000 6.2831853}%
\special{sh 1.000}%
\special{ar 1900 900 46 46  0.0000000 6.2831853}%
\special{sh 1.000}%
\special{ar 2100 900 46 46  0.0000000 6.2831853}%
\special{sh 1.000}%
\special{ar 2300 900 46 46  0.0000000 6.2831853}%

\special{ar 1600 1100 46 46  0.0000000 6.2831853}%

\special{ar 1800 1100 46 46  0.0000000 6.2831853}%

\special{ar 2000 1100 46 46  0.0000000 6.2831853}%
\special{sh 1.000}%
\special{ar 2200 1100 46 46  0.0000000 6.2831853}%
\special{sh 1.000}%
\special{ar 2400 1100 46 46  0.0000000 6.2831853}%


\special{ar 3500 300 46 46  0.0000000 6.2831853}%

\special{ar 3400 500 46 46  0.0000000 6.2831853}%
\special{sh 1.000}%
\special{ar 3600 500 46 46  0.0000000 6.2831853}%

\special{sh 1.000}%
\special{ar  3300 700 46 46  0.0000000 6.2831853}%
\special{sh 1.000}%
\special{ar 3500 700 46 46  0.0000000 6.2831853}%
\special{sh 1.000}%
\special{ar 3700 700 46 46  0.0000000 6.2831853}%

\special{ar 3200 900 46 46  0.0000000 6.2831853}%
\special{sh 1.000}%
\special{ar 3400 900 46 46  0.0000000 6.2831853}%
\special{sh 1.000}%
\special{ar 3600 900 46 46  0.0000000 6.2831853}%
\special{sh 1.000}%
\special{ar 3800 900 46 46  0.0000000 6.2831853}%

\special{ar 3100 1100 46 46  0.0000000 6.2831853}%
\special{sh 1.000}%
\special{ar 3300 1100 46 46  0.0000000 6.2831853}%
\special{sh 1.000}%
\special{ar 3500 1100 46 46  0.0000000 6.2831853}%
\special{sh 1.000}%
\special{ar 3700 1100 46 46  0.0000000 6.2831853}%
\special{sh 1.000}%
\special{ar 3900 1100 46 46  0.0000000 6.2831853}%


\special{ar 5000 300 46 46  0.0000000 6.2831853}%

\special{ar 4900 500 46 46  0.0000000 6.2831853}%
\special{sh 1.000}%
\special{ar 5100 500 46 46  0.0000000 6.2831853}%
\special{sh 1.000}%
\special{ar  4800 700 46 46  0.0000000 6.2831853}%
\special{sh 1.000}%
\special{ar 5000 700 46 46  0.0000000 6.2831853}%
\special{sh 1.000}%
\special{ar 5200 700 46 46  0.0000000 6.2831853}%

\special{sh 1.000}%
\special{ar 4700 900 46 46  0.0000000 6.2831853}%
\special{sh 1.000}%
\special{ar 4900 900 46 46  0.0000000 6.2831853}%
\special{sh 1.000}%
\special{ar 5100 900 46 46  0.0000000 6.2831853}%
\special{sh 1.000}%
\special{ar 5300 900 46 46  0.0000000 6.2831853}%

\special{ar 4600 1100 46 46  0.0000000 6.2831853}%

\special{ar 4800 1100 46 46  0.0000000 6.2831853}%
\special{sh 1.000}%
\special{ar 5000 1100 46 46  0.0000000 6.2831853}%
\special{sh 1.000}%
\special{ar 5200 1100 46 46  0.0000000 6.2831853}%
\special{sh 1.000}%
\special{ar 5400 1100 46 46  0.0000000 6.2831853}%
\end{picture}

\vspace{60pt}
\hspace{15pt}(9) \hspace{84pt} (10) 

\begin{picture}(36.47,22.63)(12.25,-30.63)

\special{ar 500 300 46 46  0.0000000 6.2831853}%
\special{sh 1.000}%
\special{ar 400 500 46 46  0.0000000 6.2831853}%

\special{ar 600 500 46 46  0.0000000 6.2831853}%
\special{sh 1.000}%
\special{ar 300 700 46 46  0.0000000 6.2831853}%
\special{sh 1.000}%
\special{ar 500 700 46 46  0.0000000 6.2831853}%
\special{sh 1.000}%
\special{ar 700 700 46 46  0.0000000 6.2831853}%
\special{sh 1.000}%
\special{ar 200 900 46 46  0.0000000 6.2831853}%
\special{sh 1.000}%
\special{ar 400 900 46 46  0.0000000 6.2831853}%
\special{sh 1.000}%
\special{ar 600 900 46 46  0.0000000 6.2831853}%
\special{sh 1.000}%
\special{ar 800 900 46 46  0.0000000 6.2831853}%

\special{ar 100 1100 46 46  0.0000000 6.2831853}%

\special{ar 300 1100 46 46  0.0000000 6.2831853}%
\special{sh 1.000}%
\special{ar 500 1100 46 46  0.0000000 6.2831853}%
\special{sh 1.000}%
\special{ar 700 1100 46 46  0.0000000 6.2831853}%
\special{sh 1.000}%
\special{ar 900 1100 46 46  0.0000000 6.2831853}%


\special{ar 2000 300 46 46  0.0000000 6.2831853}%

\special{sh 1.000}%
\special{ar 1900 500 46 46  0.0000000 6.2831853}%
\special{sh 1.000}%
\special{ar 2100 500 46 46  0.0000000 6.2831853}%

\special{sh 1.000}%
\special{ar  1800 700 46 46  0.0000000 6.2831853}%
\special{sh 1.000}%
\special{ar 2000 700 46 46  0.0000000 6.2831853}%
\special{sh 1.000}%
\special{ar 2200 700 46 46  0.0000000 6.2831853}%

\special{ar 1700 900 46 46  0.0000000 6.2831853}%
\special{sh 1.000}%
\special{ar 1900 900 46 46  0.0000000 6.2831853}%
\special{sh 1.000}%
\special{ar 2100 900 46 46  0.0000000 6.2831853}%
\special{sh 1.000}%
\special{ar 2300 900 46 46  0.0000000 6.2831853}%

\special{ar 1600 1100 46 46  0.0000000 6.2831853}%
\special{sh 1.000}%
\special{ar 1800 1100 46 46  0.0000000 6.2831853}%
\special{sh 1.000}%
\special{ar 2000 1100 46 46  0.0000000 6.2831853}%
\special{sh 1.000}%
\special{ar 2200 1100 46 46  0.0000000 6.2831853}%

\special{ar 2400 1100 46 46  0.0000000 6.2831853}%
\end{picture}

\vspace{60pt}
By Proposition \ref{prop1}, all empty circles must be at corner and each connected component of emptysets looks like the Young diaglam. 
Also, the numbers of empty circles must be less than the degree of elements of $V$.

\bigskip

\noindent
Satoshi Murai\\
Department of Pure and Applied Mathematics\\
Graduate School of Information Science and Technology\\
Osaka University\\
Toyonaka, Osaka, 560-0043, Japan\\
E-mail:s-murai@ist.osaka-u.ac.jp
\end{document}